\documentclass[12pt,a4paper]{article}
\usepackage{amssymb,amsmath,amsthm}

\newcommand\cA{{\mathcal A}}
\newcommand\cB{{\mathcal B}}
\newcommand\cC{{\mathcal C}}

\newcommand\cF{{\mathcal F}}
\newcommand\cG{{\mathcal G}}

\newcommand\cS{{\mathcal S}}

\newcommand\coF{{\overline F}}
\newcommand\coG{{\overline G}}

\newcommand\G{\Gamma}

\theoremstyle{plain}
\newtheorem{theorem}{Theorem}[section]
\newtheorem{lemma}[theorem]{Lemma}
\newtheorem{corollary}[theorem]{Corollary}
\newtheorem{conjecture}[theorem]{Conjecture}
\newtheorem{proposition}[theorem]{Proposition}
\newtheorem{problem}[theorem]{Problem}
\theoremstyle{definition}

\newcommand\lref[1]{Lemma~\ref{lem:#1}}
\newcommand\tref[1]{Theorem~\ref{thm:#1}}
\newcommand\cref[1]{Corollary~\ref{cor:#1}}

\newcommand\pref[1]{Proposition~\ref{prop:#1}}

\begin{document}

\title{Cross-Sperner families}

\author{D\'aniel Gerbner$^{{\rm a, *}}$ \and Nathan Lemons$^{{\rm b, *}}$ \and Cory
  Palmer$^{{\rm a},}$\thanks{Research 
supported by Hungarian National Scientific Fund, grant number: OTKA
    NK-78439} \and Bal\'azs Patk\'os$^{{\rm a},}$\thanks{Research 
supported by Hungarian National Scientific Fund, grant number: OTKA
    K-69062 and PD-83586} \thanks{corresponding author, e-mail: patkos@renyi.hu} \and 
Vajk Sz\'ecsi$^{{\rm b}}$
\and
\small $^{\rm a}${\it Hungarian Academy of Sciences, Alfr\'ed 
R\'enyi Institute}\\[-0.8ex]
\small {\it of Mathematics, P.O.B. 127, Budapest H-1364, Hungary}\\
\small $^{\rm b}${\it Central European University, Department of Mathematics}\\[-0.8ex]
\small {\it and its Applications, N\'ador u. 9, Budapest H-1051, Hungary}\\}

\maketitle

\begin{abstract}
A pair of families $(\cF,\cG)$ is said to be \emph{cross-Sperner} if 
there exists no pair of sets $F \in \cF, G \in \cG$ with $F \subseteq G$ 
or $G \subseteq F$. There are
two ways to measure the size of the pair $(\cF,\cG)$: with the sum
$|\cF|+|\cG|$ or with the product $|\cF|\cdot |\cG|$. We show that if $\cF,
\cG \subseteq 2^{[n]}$, then $|\cF||\cG| \le 2^{2n-4}$ and $|\cF|+|\cG|$ is 
maximal if
$\cF$ or $\cG$ consists of exactly one set of size $\lceil n/2 \rceil$ 
provided 
the size of the ground set $n$ is large enough and both $\cF$ and $\cG$ are non-empty.
\end{abstract}


\textit{keywords}: extremal set systems, Sperner property

\section{Introduction}
We use standard notation: $[n]$ denotes the set of the first $n$ positive
integers, $2^S$ denotes the power set of the set $S$ and $\binom{S}{k}$
denotes the set of all $k$-element subsets of $S$. The complement of a set $F$ is
denoted by $\coF$ and for a family $\cF$ we write $\overline{\cF}=\{\coF:F\in
\cF\}$. 

One of the first theorems in the area of extremal set families is that of
Sperner \cite{S}, stating that if we consider a family $\cF \subseteq 2^{[n]}$ 
such
that no set $F \in \cF$ can contain any other $F' \in \cF$, then the
number of sets in $\cF$ is at most ${n \choose \lfloor n/2 \rfloor}$ and
equality holds if and only if $\cF=\binom{[n]}{\lfloor n/2 \rfloor}$ or
$\cF=\binom{[n]}{\lceil n/2 \rceil}$. 
Families satisfying the
assumption of Sperner's theorem are called \emph{Sperner families} or 
\emph{antichains}.
The celebrated theorem of Erd\H os, Ko and Rado \cite{EKR} asserts that if
for a family $\cG \subseteq \binom{[n]}{k}$ we have $G \cap G' \neq 
\emptyset$ for 
all $G,G' \in \cG$ (families with this property are called 
\emph{intersecting}), 
then the size of $\cG$ is at most ${n-1 \choose k-1}$ provided $2k \le n$.

There have been many generalizations and extensions both to the 
theorem of Sperner and to the result by 
Erd\H os, Ko and Rado (two excellent but not really recent surveys are 
\cite{DF} and \cite{E}). One such generalization is the following: a 
pair $(\cF,\cG)$ of
families is said to be \emph{cross-intersecting} if for any $F \in \cF, G 
\in \cG$ we have 
$F \cap G\neq \emptyset$. Cross-intersecting pairs of families have been 
investigated for quite a while
and attracted the attention of many researchers 
\cite{Be,B,FT,FT2,Fu,K1,K2,KS}. The present paper deals with 
the analogous 
generalization of Sperner families that has not been considered in the 
literature. A pair $(\cF,\cG)$ of
families is said to be \textit{cross-Sperner} if 
there exists no pair of sets $F \in \cF, G \in \cG$ with $F \subseteq G$ 
or $G \subseteq F$. There are two ways to measure the size of the 
pair $(\cF,\cG)$: either with the sum
$|\cF|+|\cG|$ or with the product $|\cF|\cdot |\cG|$. We will address both
problems. 

Clearly, $|\cF|+|\cG| \le 2^n$ as by definition $\cF \cap \cG=\emptyset$. 
 The sum $2^n$ can be obtained  by putting $\cF=\emptyset, \cG=2^{[n]}$. 
Thus, when considering 
the problem of maximizing $|\cF|+|\cG|$ we will assume that both $\cF$ and 
$\cG$ are non-empty.

We can reformulate our problem in a rather interesting way. Let 
$\G_n=(V_n,E_n)$ be the graph with vertex set $V_n=2^{[n]}$
and edge set $E_n=\{(F,G):F,G \in V_n ~\ F \subsetneq G$ or 
$G \subsetneq F\}$. Then $\max\{|\cF|+|\cG|\}=2^n-c(\G_n)$, where
$c(\G_n)$ denotes the vertex connectivity of $\G_n$. Moreover, if we let
\[
F(n,m)=\max\{|\cG|:\cG \subseteq 2^{[n]}, \exists \cF \subseteq 
2^{[n]} ~\text{with}\ |\cF|=m, (\cF,\cG) ~\text{is cross-Sperner}\},
\] 
then, denoting by $N_{\G_n}(U)$ the neighborhood of $U$ in $\G_n$, we 
have 
\[
  F(n,m)=2^n-m-\min\{|N_{\G_n}(\cF)|: \cF \subseteq V_n, |\cF|=m\}.
\] 
Thus determining
$F(n,m)$ is equivalent to the isoperimetric problem for the graph $\G_n$.

Let us mention that the cross-Sperner property of the pair $(\cF,\cG)$ is
equivalent to $(\cF, \overline{\cG})$ being cross-intersecting and
cross-co-intersecting, i.e. for any $F \in \cF$ and $G \in \cG$ we have $F
\cap \coG \neq \emptyset$ and $F \cup \coG\neq [n]$.

The rest of the paper is organized as follows. In Section 2, we consider the
problem of maximizing $|\cF|+|\cG|$ and prove the following theorem.
\begin{theorem}
\label{thm:sum} There exists an integer $n_0$ such that if $n \ge n_0$ and
the pair $(\cF,\cG)$ is cross-Sperner with $\emptyset \neq \cF,\cG 
\subseteq 2^{[n]}$, then 
\[
|\cF|+|\cG| \le F(n,1)+1=2^n-2^{\lceil n/2\rceil}-2^{\lfloor n/2\rfloor}+2,
\]
and equality holds if and only if $\cF$ or $\cG$ consists of exactly one set 
$S$ of size $\lfloor n/2 \rfloor$ or $\lceil n/2 \rceil$  and the other family
consists of all subsets of $[n]$ not contained in $S$ and not containing $S$.
\end{theorem}
In Section 3, we address the problem of maximizing $|\cF|\cdot |\cG|$. Our
result is the following theorem.
\begin{theorem}
\label{thm:prod} If $n \ge 2$ and $(\cF,\cG)$ is cross-Sperner with $\cF,\cG 
\subseteq 2^{[n]}$, then the 
following inequality holds:
$$|\cF||\cG|\le 2^{2n-4}.$$
This bound is best possible as shown by $\cF=\{F\in 2^{[n]}: 1 \in 
F, n \notin F\},
\cG=\{G\in 2^{[n]}: n \in G, 1 \notin G\}$.
\end{theorem}
Finally, Section 4 contains some concluding remarks and open problems.

\section{Proof of \tref{sum}}

Before we start the proof of \tref{sum}, let us introduce some notation and
state a theorem that we will use in our proof. For a $k$-uniform family $\cF
\subseteq \binom{[n]}{k}$ let $\Delta \cF=\{G \in \binom{[n]}{k-1}: \exists F
\in \cF, G \subset F\}$ be the \emph{shadow} of $\cF$. The following version
of the shadow theorem
is due to Lov\'asz \cite{L}.

\begin{theorem}
\label{thm:lovasz} \emph{[Lov\'asz \cite{L}]} Let $\cF
\subseteq \binom{[n]}{k}$ and let us define the real number $x$ by 
$|\cF|=\binom{x}{k}$. Then we have $\Delta \cF \ge \binom{x}{k-1}$.
\end{theorem}

For any $F \in 2^{[n]}$ we have $N_{\G_n}(F)=2^{|F|}+2^{n-|F|}-2$ which is
minimized if $|F|=\lceil n/2 \rceil$. This proves $ F(n,1)=2^n-2^{\lceil
  n/2\rceil}-2^{\lfloor n/2\rfloor}+1$ as stated in \tref{sum}.

\begin{proposition}
\label{prop:convex} If a pair $(\cF,\cG)$ maximizes $|\cF|+|\cG|$, then both
$\cF$ and $\cG$ are convex families i.e. $F_1 \subset F \subset F_2$,
$F_1,F_2 \in \cF$ implies $F \in \cF$.
\end{proposition}

\begin{proof}
If $F,F_1,F_2$ are as above, then $F$ can be added to $\cF$ since any set
containing $F$ contains $F_1$ and any subset of $F$ is a subset of
$F_2$. 
\end{proof}

Let $(\cF,\cG)$ be a pair of cross-Sperner families and let $F_0$ and
$G_0$ be sets of minimum size in $\cF$ and $\cG$.

\begin{proposition}
\label{prop:small} If $|F_0|+|G_0| < \lceil n/2 \rceil -1$, then $|\cF|+|\cG|
< F(n,1)$.
\end{proposition}

\begin{proof}
No set containing $F_0 \cup G_0$ can be a member of $\cF$ or $\cG$.
\end{proof}

As $(\cF,\cG)$ is cross-Sperner if and only if
$(\overline{\cF},\overline{\cG})$ is cross-Sperner, by taking complements (if
necessary) and \pref{small} we may and will assume that $m:=|F_0|\ge \lfloor 
n/4 \rfloor$. Let us write $\cF^*=\{F \in \cF: F_0 \subsetneq F\}$. Subsets of
$F_0$ are not in $\cF$ by the minimality of $F_0$ and by the cross-Sperner
property  they 
cannot be in $\cG$ either, thus to prove \tref{sum} we need to show that there
exist more than $|\cF^*|$ many sets that are not contained in $\cF \cup \cG$
and are not subsets of $F_0$. For any $F^* \in \cF^*$ let us define 
\[
 B(F^*)=\{F^*\setminus F'_0: F'_0 \subseteq F_0, |F^*\setminus F'_0|<m\}.
\]
Clearly, for any $F^*_1,F^*_2 \in \cF^*$ we have $B(F^*_1) \cap
B(F^*_2)=\emptyset$ as they already differ outside $F_0$. By definition, no
set in $\cB:=\cup_{F^* \in \cF^*}B(F^*)$ is a subset of $F_0$. We have $\cB
\cap \cF=\emptyset$ as all sets in $\cB$ have size smaller than $m$ and $\cB
\cap \cG=\emptyset$ by the cross-Sperner property. Thus to prove \tref{sum} it
is enough to show that $|\cF^*| < |\cB|$.

Note the following three things:
\begin{itemize}
\item
$|B(F^*)|=\sum_{i=|F^*\setminus F_0|+1}^m\binom{m}{i}$,
\item
$\cF^{**}=\{F^*\setminus F_0:F^* \in \cF^*\}$ is downward closed as $\cF$ and
$\cF^*$ are convex,
\item
$|\cF^{**}|=|\cF^*|$.
\end{itemize}

Therefore the following lemma finishes the proof of \tref{sum} by choosing
$\cA=\cF^{**}$, $k=m$ and $n'=n-|F_0|$.

\begin{lemma} 
\label{lem:sterling}Let $\emptyset \ne \cA\subseteq 2^{[n']}$ be a downward 
closed family and 
$k\ge n'/3$. Then if $n'$ is large enough, the following holds

\begin{equation} 
\label{1}|\cA|< \sum_{A\in\cA} \sum_{i=|A|+1}^k {k\choose i}. 
\end{equation}

\end{lemma}

\begin{proof} Let $a_i=|\{A\in\cA: |A|=i\}|$ and $w(j)=\sum_{i=j+1}^k 
{k\choose i}$. Then we can formulate (\ref{1}) in the following way:

\begin{equation} \sum_{j=0}^{n'} a_j< \sum_{j=0}^{n'} a_jw(j). \end{equation}

Let $x$ be defined by $a_{k-1}={x \choose k-1}$. By \tref{lovasz} 
if $j< k-1$ 
then $a_j\ge {x\choose j}$. If we replace $a_j$
by ${x\choose j}$ in (2), then the LHS decreases by $a_j-{x\choose
  j}$ and the RHS decreases by $(a_j-{x\choose j})w(j)$, which is
larger. If $j> k-1$, then $a_j\le {x\choose j}$ again by \tref{lovasz}. 
If we replace $a_j$ by ${x\choose j}$ in (2), then the
LHS increases while the RHS does not change (as for $j \ge
k$ we have $w(j)=0$). 
Hence it is enough to prove 

\begin{equation}
\label{3} \sum_{j=0}^{n'} {x\choose j}< \sum_{j=0}^{n'} {x\choose j}w(j). 
\end{equation}

First we prove (\ref{3}) for $x=n'$. In this case the LHS is $2^{n'}$ while 
the RHS is
monotone increasing in $k$, thus it is enough to prove for $k=\lceil n/3\rceil$. We will 
estimate the RHS from below by considering only one term of the sum. Clearly, 
$\binom{n'}{j}w(j) \ge \binom{n'}{j}\binom{k}{j+1} \ge \binom{n'}{j}\binom{n'/3}{j+1}$. 
Let us write $j=\alpha n'$ for some $0\le \alpha \le 1/3$. Then by Stirling's formula we obtain 
\[
\binom{n'}{j}\binom{n'/3}{j+1}=\binom{n'}{\alpha n'}\binom{n'/3}{\alpha n'+1}=\Theta\left(\frac{1}{n'}\left(\frac{1}{\alpha^{2\alpha}(1-\alpha)^{1-\alpha}3^{1/3}(1/3-\alpha)^{1/3-\alpha}}\right)^{n'}\right).
\]
The value of the fraction in parenthesis is larger than 2 for, say, $\alpha=2/9$, thus (\ref{3}) holds if $n'$ is large enough and $x=n'$.

\vskip 0.2truecm

To prove (\ref{3}) for arbitrary $x$, let $c={x\choose k-1}/{n'\choose k-1}$. If $j>k-1$, 
then $c>{x\choose j}/{n' \choose j}$, while if $j<k-1$, then 
$c<{x\choose j}/{n' \choose j}$. By the $x=n'$ case we know 

\begin{equation} \sum_{j=0}^{n'} c{n'\choose j}< \sum_{j=0}^{n'} c{n'\choose
    j}w(j). 
\end{equation}

Let us replace $c{n'\choose j}$ by ${x\choose j}$ in this inequality. 
If $j>k-1$, then the LHS decreases and the RHS does not 
change. If $j=k-1$ none of the sides change by definition of $c$. If $j<k-1$, 
both sides increase, and 
the RHS increases more as $w(j) \ge 1$ for all $0 \le j \le k-1$. 
Hence the inequality holds and 
gives back (3), which finishes the proof of the lemma.
\end{proof}

We believe that \tref{sum} is valid for all $n$, but unfortunately
\lref{sterling} fails for small values of $n$.

\section{Proof of \tref{prod}}

In this section we prove \tref{prod}. Our main tool will be the following
special case of the Four Functions Theorem of Ahlswede and Daykin \cite{AD}. To state
their result for any pair $\cA,\cB$ of families let us write $\cA \wedge 
\cB=\{A \cap B: A \in \cA, B \in \cB\}$ and
$\cA \vee \cB=\{A \cup B: A \in \cA, B \in \cB\}$.

\begin{theorem}
\label{thm:FFT}\emph{[Ahlswede-Daykin, 
\cite{AD}]}
For any pair $\cA,\cB$ of families  we have 
$$|\cA||\cB| \le |\cA \wedge \cB||\cA \vee \cB|.$$
\end{theorem}

To prove \tref{prod} we will need the following lemma.

\begin{lemma}
\label{lem:partition} If $(\cF,\cG)$ is a pair of cross-Sperner families, then
the families $\cF$, $\cG$, $\cF \wedge \cG$ and $\cF \vee \cG$ are 
pairwise disjoint.
\end{lemma}

\begin{proof}
$\cF$ and $\cG$ are disjoint as some set $F \in \cF \cap \cG$ is a subset of
itself and thus contradicts the
cross-Sperner property.
$\cF$ and $\cG$ are both disjoint from $\cF \wedge \cG$ and $\cF \vee \cG$ as 
$F \cap G \subseteq F,G$ and $F,G \subseteq F \cup G$.
Finally, $\cF \wedge \cG$ and $\cF \vee \cG$ are disjoint as 
$F_1 \cap G_1=F_2 \cup G_2$ would imply $F_2 \subseteq G_1$.
\end{proof}

Now we are able to prove \tref{prod}. 

\begin{proof}
Let $(\cF,\cG)$ be a cross-Sperner pair
of families. Clearly, if $|\cF|+|\cG| \le 2^{n-1}$, then the statement 
of the theorem holds.
But if $|\cF|+|\cG| > 2^{n-1}$, then by \lref{partition} we have 
$|\cF \wedge \cG|+|\cF \vee \cG|<2^{n-1}$ and
thus by \tref{FFT} we obtain $|\cF||\cG| \le |\cF \wedge \cG||\cF \vee \cG| 
\le 2^{2n-4}$.
\end{proof}

\begin{corollary}
\label{cor:2n2} For $n \ge 2$, we have $F(n,2^{n-2})=2^{n-2}$.
\end{corollary}

%With a little extra work one can prove that the same bound on $|\cG||\cG|$ 
%holds if one allows sets to be members of both $\cF$ and $\cG$.

%\begin{lemma}
%\label{lem:spernershadow}
%Let $\cA \subseteq 2^{[n]}$ be an antichain and let us define $x$ by $\lceil
%|\cA|/2 \rceil = \binom{x}{\lceil n/2 \rceil}$. Then the following inequality
%holds: 
%\[
% |N(\cA)| \ge \sum_{i=0}^{\lceil n/2 \rceil-1}\binom{x}{i}.
%\]
%\end{lemma}

%\begin{proof}
%W.l.o.g we may assume that $\cA^+:=\cA \cap \binom{[n]}{\ge \lceil 
%n/2 \rceil}$ has size at least $\binom{x}{\lceil n/2 \rceil}$. Clearly, we
%have $|\Delta_{\lceil n/2 \rceil}(\cA^+)|\ge |\cA^+|$ and then we are done by
%\tref{}. 
%\end{proof}

%\begin{corollary}
%\label{cor:spernershadow2}
%Let $\cA \subseteq 2^{[n]}$ be an antichain. If $n \ge ...$, then we have
%\[
%  |N(\cA)| \ge 3|\cA|.
%\]
%\end{corollary}

%\begin{proof}
%Let $x$ be as in \lref{spernershadow}. Then if $x \le n-12$, then we are done
%as in that case $\binom{x}{\lceil n/2 \rceil} \le \binom{x}{\lceil n/2 
%\rceil-i}$ for $1\le i\le 6$. While if $x \ge n-12$, then by
%\lref{spernershadow} we obtain
%\[|N(\cA)|\ge 2^{n-13} \ge 3\binom{n}{\lceil n/2 \rceil},
%\]
%provided $n$ is large enough. Then we are done as by Sperner's theorem $|\cA|
%\le \binom{n}{\lceil n/2 \rceil}$.
%\end{proof}

%Now we are ready to prove the analog of \tref{prod} for cross-Sperner 
%pairs $(\cG,\cF)$ of families where $\cF \cap \cG \neq \emptyset$ is allowed.

%\begin{theorem}
%\label{thm:prod+}
%\end{theorem}

\section{Concluding remarks and open problems}

One might wonder whether it changes the situation
if we allow sets to belong to both $\cF$ and $\cG$ and we modify the 
definition of
cross-Sperner families so that only pairs $F \in \cF, G \in \cG$ with 
$F \subsetneq G$ or
$G \subsetneq F$ are forbidden. It is easy to see that the situation is the
 same when
considering $|\cF|+|\cG|$. To prove that $|\cF|+|\cG| \le 2^n$ let us write 
$\cC=\cF \cap \cG$ and if it is not empty,
then $D(\cC):=
\{C\setminus C': C,C' \in \cC\}$ is disjoint both from $\cF$ and $\cG$ and 
a result by Marica and Sch\"onheim \cite{MS}
tells us that $|D(\cC)| \ge |\cC|$. Note that the proof of \tref{sum} works in
this case as well giving the upper bound $|\cF|+|\cG| \le F(n,1)+2$.

Although $F(n,m)$ is not known for most values, it is natural to generalize 
the problem to $k$-tuples of
families: $\cF_1,\cF_2,...,\cF_k$ is said to be cross-Sperner if for any 
$1\le i<j\le k$ there is no pair $F \in \cF_i$ and $F' \in \cF_j$ with 
$F \subseteq F'$ or $F' \subseteq F$. One can consider the problems of 
maximizing $\sum_{i=1}^k|\cF_i|$ and
$\prod_{i=1}^k|\cF_i|$. In the former case we need the extra assumption that
all $\cF_i$ are non-empty as otherwise the trivial upper bound $2^n$ is tight.

When maximizing the sum, it is natural to conjecture that in the best possible 
construction all but one family consists of one single set. By the 
cross-Sperner property, these sets together must form a Sperner family, 
therefore it might turn out to be useful to introduce
\[
F^*(n,m)=\max\{|\cG|:\cG \subseteq 2^{[n]}, \exists \cF \subseteq 
2^{[n]} \]
\[~\text{with}\ |\cF|=m, (\cF,\cG) ~\text{is cross-Sperner, $\cF$ is Sperner}\}.
\] 
\begin{problem}
Under what conditions is it true that if $\cF_1,\cF_2,...,\cF_k$ form a $k$-tuple of non-empty cross-Sperner families, then
$$\sum_{i=1}^k|\cF_i| \le k-1+F^*(n,k-1)?$$
\end{problem}
Concerning maximizing the product of the $|\cF_i|$, by \tref{prod} one 
obtains that 
\[
\prod_{i=1}^k|\cF_i|=\left(\prod_{1\le i<j\le k}|\cF_i||\cF_j|\right)^{\frac{1}{k-1}}\le 2^{kn-2k}.
\]
We conjecture that the following construction is optimal: let $l=l(k)$ be 
the smallest positive integer so that $k \le \binom{l}{\lfloor l/2 \rfloor}$. 
Then there exists a Sperner family $\cS=\{S_1,...,S_k\}\subseteq 2^{[l]}$ of 
size $k$. Put
$\cF_i=\{F \subseteq [n]: F \cap [l]=S_i\}$. Clearly, the $\cF_i$ form a 
$k$-tuple of cross-Sperner families and we have 
$\prod_{i=1}^k|\cF_i|=2^{k(n-l)}$. Unfortunately, already for $l=3$ there is 
a gap of a factor of 8 between the upper bound and the size of our 
construction.
\begin{conjecture}
If $\cF_1,\cF_2,...,\cF_k\subseteq 2^{[n]}$ form a $k$-tuple of cross-Sperner families, then
\[
\prod_{i=1}^k|\cF_i| \le 2^{k(n-l)},
\]
where $l$ is the least positive integer with $\binom{l}{\lfloor l/2 \rfloor} \ge k$.
\end{conjecture}

\end{document}